\documentclass[a4paper,12pt]{article}
\usepackage{amsmath,amsthm}
\usepackage[letterpaper,margin=1in]{geometry}
\usepackage{blindtext}
\usepackage[symbol]{footmisc}
\numberwithin{equation}{section}

\usepackage[colorlinks,citecolor=blue,urlcolor=red,bookmarks=false,hypertexnames=true]{hyperref}
\DeclareRobustCommand{\rchi}{{\mathpalette\irchi\relax}}
\newcommand{\irchi}[2]{\raisebox{\depth}{$#1\chi$}} 
\begin{document}
\centerline{}\centerline{\Large{\bf A Study On Some Geometric Properties and Physical}}
\centerline{}\centerline{\Large{\bf Applications of A Mixed Quasi-Einstein Spacetime }}
\centerline{}\centerline{$~$Kaushik Chattopadhyay\textsuperscript{1}\footnote{This is the corresponding author.}, Arindam Bhattacharyya\textsuperscript{2} and Dipankar Debnath\textsuperscript{3}}
\centerline{\textsuperscript{1}Department of Mathematics,}\centerline{Jadavpur University, Kolkata-700032, India}
\centerline{E-mail:kausikchatterji@gmail.com}
\centerline{\textsuperscript{2}Department of Mathematics,}\centerline{Jadavpur University, Kolkata-700032, India}
\centerline{E-mail:bhattachar1968@yahoo.co.in}
\centerline{\textsuperscript{3}Department of Mathematics,}\centerline{Bamanpukur High School(H.S), Nabadwip, India}
\centerline{E-mail:dipankardebnath123@hotmail.com}
\newtheorem{Theorem}{Theorem}[section]
\begin{abstract}{\em In the present paper we discuss about a set of geometric properties and physical applications of a mixed quasi-Einstein spacetime$[M(QE)_{4}]$, which is a special type of nearly quasi-Einstein spacetime$[N(QE)_{4}]$. Firstly we consider a nearly quasi-Einstein manifold$[N(QE)_{n}]$ along with a mixed quasi-Einstein manifold$[M(QE)_{n}]$ and study some geometric conditions with respect to different curvature tensors on them. Then we consider a mixed quasi-Einstein spacetime and study $W_{2}$-Ricci pseudosymmetry and projective Ricci pseudosymmetry on it. Then we study the conditions $H \cdot S=0$ and $\tilde C \cdot S=0$ in an $M(QE)_{4}$ spacetime, where $H, \tilde C$ are the conharmonical curvature tensor and semiconformal curvature tensor respectively. Next we consider a semiconformally flat Ricci pseudosymmetric $M(QE)_{4}$ spacetime and find the nature of that spacetime. In the next section we study about some applications of a perfect fluid $M(QE)_{4}$ spacetime in the general theory of relativity. Followed by we discuss about the Ricci soliton structure in a semiconformally flat $M(QE)_{4}$ spacetime satisfying Einstein field equation without cosmological constant. Finally we give a nontrivial example of a $M(QE)_{4}$ spacetime to establish the existence of it.}
\end{abstract}
\textbf{M.S.C.2010:} ~~53C15, 53C25, 53C35.

\section{Introduction }

$~~~~$The General theory of Relativity is unarguably the most beautiful theory the world of Physics has ever produced. It is the most powerful result of the human intellect. This is an extremely important theory to study the nature of this universe, cosmology and gravity. A crucial thing that Modern Scientists and Mathematical Physicists can learn from special to general relativity is that, the spacetime has to be considered as a class of semi-Riemannian geometry.
Thus the semi-Riemannian geometry has become more and more relevant and significant in dealing with the nature of this universe with every passing day.\\\\
$~~~~$In general relativity the matter content of the spacetime is described by the energy momentum tensor $T$. Since the matter content of the universe is assumed to behave like a perfect fluid in the standard cosmological models, the physical motivation for studying Lorentzian manifolds is the assumption that a gravitational field may be effectively modelled by some Lorentzian metric defined on a suitable four dimensional manifold $M$. The study of different spacetimes arose while describing the different phases of the evolution of the universe. After the big bang our universe has evolved through $3$ different phases of evolution. Namely,\\\\
\textbf{\em i) The initial phase:}  This is just after the big bang when both the viscosity and the heat flux of the universe were present.\\
\textbf{\em ii) The intermediate phase:} In this stage the viscosity became negligible but the heat flux was still non-negligible.\\
\textbf{\em iii) The final phase:} In this stage both the viscosity and the heat flux have become negligible. In the present days the universe belongs to this phase.\\\\
Quasi-Einstein spacetime($\cite{jmh2}$) describes the basic structure of the final phase and generalized quasi-Einstein spacetime($\cite{jmh56}$, $\cite{jmh55}$) describes the structure of the universe in the intermediate phase.\\\\
$~~~~$The mathematical structure of this universe is mainly studied on a semi-Riemannian manifold which sometimes is not an Einstein spacetime $\cite{jmh5}$. Like the Robertson-Walker spacetime is not an Einstein spacetime $\cite{jmh52}$, rather it's a quasi-Einstein manifold $\cite{jmh53}$. Thus it was always necessary to extend the concept of Einstein manifold to quasi-Einstein. But a manifold may not be even quasi-Einstein, yet nearly equal to a quasi-Einstein manifold. Those are named as nearly quasi-Eintstein manifolds $\cite{jmh44}$. To be more specific an Einstein manifold is always a quasi-Einstein manifold but the converse is not true in general. And a quasi-Einstein manifold is always a nearly quasi-Einstein manifold but the converse is not true in general. The mixed quasi-Einstein spacetime works as a bridge between an Einstein spacetime and a generalized quasi-Einstein spacetime ($\cite{jmh56}$, $\cite{jmh3}$) which inspired us to take a look on mixed quasi-Einstein spacetime.\\\\
$~~~~$The present paper is organised as follows:\\\\
After the preliminaries in section (\ref{s2}) we derive the natures of the $J$-flat nearly quasi-Einstein manifolds and $J$-flat mixed quasi-Einstein manifolds(where $J$ is a type of curvature tensor) in section (\ref{s3}). Then we discuss about some geometric properties of $M(QE)_{4}$ spacetime, namely $W_{2}$-Ricci pseudosymmetry in section (\ref{s4}) and projective Ricci pseudosymmetry in section (\ref{s5}). Next we study the conditions $H \cdot S=0$ and $\tilde C \cdot S=0$ in an $M(QE)_{4}$ spacetime in section (\ref{s6}). After that we discuss about the nature of the semiconformally flat Ricci pseudosymmetric $M(QE)_{4}$ spacetime in section (\ref{s7}). Next we discuss about some applications of a perfect fluid $M(QE)_{4}$ spacetime in general relativity in section (\ref{s8}). Then we discuss about the Ricci soliton structure of a semiconformally flat $M(QE)_{4}$ spacetime in section (\ref{s9}). Finally in section (\ref{s10}) we construct a nontrivial example of an $M(QE)_{4}$ spacetime to establish the existence of it.
\section{Preliminaries}\label{s2}
$~~~$An Einstein manifold is a Riemannian or pseudo-Riemannian manifold whose Ricci tensor S of type $(0, 2)$ is non-zero and proportional to the metric tensor. Einstein manifolds form a natural subclass of various classes of Riemannian or semi-Riemannian manifolds by a curvature condition imposed on their Ricci tensor $\cite{jmh1}$. Also in Riemannian geometry as well as in general relativity theory, the Einstein manifold plays a very important role.\\\\
$~~~~$M. C. Chaki and R. K. Maity had given the notion of quasi-Einstein manifold $\cite{jmh2}$ in $2000$. A non flat $n$-dimensional Riemannian manifold $(M^{n},g)$, $n(>2)$ is said to be a quasi-Einstein manifold if its nonzero Ricci tensor $S$ of type $(0,2)$ satisfies the following condition
\begin{equation} \label{2.1}
S(X,Y)=\alpha g(X,Y)+\beta A(X)A(Y),
\end{equation}
where for all vector fields $X$,
\begin{equation} \label{2.2}
g(X,\xi_{1})= A(X),~g(\xi_{1},\xi_{1})= 1.
\end{equation}
That is, $A$ being the associated $1$-form, $\xi_{1}$ is generally known as the generator of the manifold. $\alpha$ and $\beta$ are associated scalar functions such that $\beta \neq 0$. This manifold is denoted by  $(QE)_n$. Clearly, for $\beta= 0$, this manifold reduces an Einstein manifold. We can note that Robertson-Walker spacetimes are quasi-Einstein spacetimes. In the paper $\cite{jmh4}$ A. A. Shaikh, D. Y. Yoon and S. K. Hui have studied about the properties and the application of quasi-Einstein spacetime. Many more works have been done in the spacetime of  general relativity $\cite{jmh31}$, $\cite{jmh41}$, $\cite{jmh62}$, $\cite{jmh22}$, $\cite{jmh3}$ $\cite{jmh7}$, $\cite{jmh8}$, $\cite{jmh38}$, $\cite{jmh78}$, $\cite{jmh27}$, $\cite{jmh28}$, $\cite{jmh5}$.\\\\
$~~~~$In the year of $2008$ A. K. Gazi and U. C. De $\cite{jmh44}$ had defined a new type of manifold called a nearly quasi-Einstein manifold. The very next year in another paper $\cite{jmh45}$ they have established the existence theorem of such type of manifold. According to them a non-flat semi-Riemannian manifold $M_{n}(n \ge 3)$ is said to be a nearly quasi-Einstein manifold if its Ricci tensor satisfies the following equation:\\
\begin{equation}\label{2.3}
S(X,Y)=\alpha g(X,Y)+\beta D(X,Y),
\end{equation}
where $\alpha, \beta$ are real-valued scalar functions such that $\beta \neq 0$ and $D$ is a nonzero symmetric tensor of type $(0,2)$. A nearly quasi-Einstein manifold is generally denoted as $N(QE)_{n}$. There are many works have been done on this type of manifold like $\cite{jmh47}$, $\cite{jmh46}$, $\cite{jmh43}$ and $\cite{jmh48}$.\\\\
$~~~~$In $\cite{jmh46}$ R. N. Singh, M. K. Pandey and D. Gautam have defined a type of nearly quasi-Einstein manifold by choosing the tensor $D$ as follows:\\
\begin{equation}\label{2.4}
D(X,Y)=A(X)B(Y)+A(Y)B(X),
\end{equation}
where $g(X,\xi_{1})= A(X)$,~$g(X,\xi_{2})= B(X)$,~$g(\xi_{1},\xi_{1})= 1$,~$g(\xi_{2},\xi_{2})= 1$,~$g(\xi_{1},\xi_{2})= 0$, for all $X \in \rchi(M)$. Combining the equations $(\ref{2.3})$ and $(\ref{2.4})$ we get
\begin{equation}\label{2.5}
S(X,Y)=\alpha g(X,Y)+\beta [A(X)B(Y)+A(Y)B(X)],
\end{equation}
where $\alpha, \beta$ are real-valued scalar functions such that $\beta \neq 0$ and $g(X,\xi_{1})= A(X)$,~$g(X,\xi_{2})\\= B(X)$,~$g(\xi_{1},\xi_{1})= 1$,~$g(\xi_{2},\xi_{2})= 1$,~$g(\xi_{1},\xi_{2})= 0$, for all $X \in \rchi(M)$. Taking a contraction over $X, Y$ from (\ref{2.5}) we get
\begin{equation}\label{2.6}
r=4\alpha,
\end{equation}
where $r$ is called the scalar curvature of the manifold.\\\\
$~~~~$In $2010$ H. G. Nagaraja had introduced the concept of a mixed quasi-Einstein manifold $\cite{jmh63}$. Later S. Mallick and U. C. De ($\cite{jmh49}$, $\cite{jmh71}$) had modified the definition of a mixed quasi-Einstein manifold given by Nagaraja and they have defined the same manifold satisfying $(\ref{2.5})$ as a mixed quasi-Einstein manifold[$M(QE)_{n}$]. They have also studied about some properties of a mixed quasi-Einstein spacetime $\cite{jmh69}$, $\cite{jmh38}$. Thus we can draw the following remark:\\\\
\textbf{Remark 2.1:} {\em A mixed quasi-Einstein manifold is a special type of nearly quasi-Einstein manifold.}\\\\
In case of we consider a mixed quasi-Einstein spacetime, it will be denoted as $M(QE)_{4}$ spacetime and the conditions will be reduced to
\begin{equation}\label{2.7}
g(X,\xi_{1})= A(X),g(X,\xi_{2})= B(X),g(\xi_{1},\xi_{1})=-1,g(\xi_{2},\xi_{2})=1,g(\xi_{1},\xi_{2})= 0,
\end{equation}
for all $X \in \rchi(M)$.\\\\
$~~~~$In $2009$ U. C. De and A. K. Gazi $\cite{jmh45}$ have introduced the notion of a manifold of nearly quasi-constant curvature. According to them, a Riemannian manifold is called a manifold of nearly quasi-constant curvature, if the curvature tensor $\tilde R$ of type $(0,4)$ satisfies the following condition:
\begin{eqnarray}\label{2.8}
&&\tilde R(X,Y,Z,W)=a[g(Y,Z)g(X,W)-g(X,Z)g(Y,W)]\nonumber\\
&&~~~~~~~~~~~~~~~~~~~~+b[g(X,W)E(Y,Z)-g(X,Z)E(Y,W)\nonumber\\
&&~~~~~~~~~~~~~~~~~~~~+g(Y,Z)E(X,W)-g(Y,W)E(X,Z)],
\end{eqnarray}
where $\tilde R(X,Y,Z,W)=g(R(X,Y)Z,W),~\tilde R$ is the curvature tensor of type $(0,4)$, $a,~b$ are scalar functions and $E$ is a nonzero symmetric tensor of type $(0,2)$. They denoted such a manifold as $N(QC)_{n}$.\\\\
$~~~~$A non-flat Riemannian manifold $(M^{n},g),(n\ge3)$ is called a manifold of mixed quasi-constant curvature $\cite{jmh71}$ if its curvature tensor $\tilde R$ of type $(0, 4)$ satisfies the condition:
\begin{eqnarray}\label{2.9}
&&\tilde R(X,Y,Z,W)=a[g(Y,Z)g(X,W)-g(X,Z)g(Y,W)]\nonumber\\
&&~~~~~~~~~~~~~~~~~~~~+b[g(X,W)\{A(Y)B(Z)+A(Z)B(Y)\}\nonumber\\
&&~~~~~~~~~~~~~~~~~~~~-g(X,Z)\{A(Y)B(W)+A(W)B(Y)\}\nonumber\\
&&~~~~~~~~~~~~~~~~~~~~+g(Y,Z)\{A(X)B(W)+A(W)B(X)\}\nonumber\\
&&~~~~~~~~~~~~~~~~~~~~-g(Y,W)\{A(X)B(Z)+A(Z)B(X)\}],
\end{eqnarray}
where $\tilde R(X,Y,Z,W)=g(R(X,Y)Z,W),~\tilde R$ is the curvature tensor of type $(0,4)$, $a,~b$ are scalar functions, $A,B$ are nonzero $1$-forms. A manifold of mixed quasi-constant curvature is denoted as $M(QC)_{n}$. If we consider a nonzero symmetric tensor of type $(0,2)$, $E'$ such that $E'(X,Y)=A(X)B(Y)+A(Y)B(X)$ then equation (\ref{2.9}) can be written as:
\begin{eqnarray}\label{2.10}
&&\tilde R(X,Y,Z,W)=a[g(Y,Z)g(X,W)-g(X,Z)g(Y,W)]\nonumber\\
&&~~~~~~~~~~~~~~~~~~~~+b[g(X,W)E'(Y,Z)-g(X,Z)E'(Y,W)\nonumber\\
&&~~~~~~~~~~~~~~~~~~~~+g(Y,Z)E'(X,W)-g(Y,W)E'(X,Z)].
\end{eqnarray}
From (\ref{2.8}), (\ref{2.9}) and (\ref{2.10}) we draw the following remark,\\\\
\textbf{Remark 2.2} {\em An $M(QC)_{n}$ is always an $N(QC)_{n}$.}\\\\
$~~~~$Now we give the expressions of the 4 important curvature tensors which we have rigorously used in our research work.\\\\
$~~~~$The $W_{2}$-curvature tensor was introduced by G. P. Pokhariyal and R. S. Mishra $\cite{jmh10}$ and they studied some properties of it. The $(0,4)$ $W_{2}$-curvature tensor on a manifold  $M^{n}(n>3)$ is defined by
\begin{equation} \label{2.11}
W_{2}(X,Y,Z,W)=\tilde R(X,Y,Z,W)-\frac{1}{n-1}[g(Y,Z)S(X,W)-g(X,Z)S(Y,W)],
\end{equation}
$\forall~X, Y, Z, W \in \rchi(M).$ Here $\tilde R$ being the curvature tensor and $Q$ is the Ricci operator defined by $g(QX,Y)= S(X,Y),\forall~X,Y \in \rchi(M)$.\\\\
$~~~~$The projective curvature tensor plays a very significant role in studying different properties of semi-Riemannian geometry. Let $M^{n}(n \ge 3)$ be a semi-Riemannian manifold. The $(0,4)$ projective curvature tensor$\cite{jmh32}$ is defined by
\begin{equation}\label{2.12}
P(X,Y,Z,W)=\tilde R(X,Y,Z,W)-\frac{1}{n-1}\{S(Y,Z)g(X,W)-S(X,Z)g(Y,W)\},
\end{equation}
$\forall~X, Y, Z, W \in \rchi(M).$ and $\tilde R$ is the $(0,4)$ Riemannian curvature tensor.\\\\
$~~~~$The conformal curvature tensor on a manifold  $M^{n}(n>3)$ was introduced by H. Weyl($\cite{jmh50}$, $\cite{jmh54}$). It was found by using the tool of conformal transformation. The tensor vanishes whenever the metric is conformally flat(for this reason the tensor was named as conformal curvature tensor). It is given by
\begin{eqnarray}\label{2.13}
&&C(X,Y,Z,W)=\tilde R(X,Y,Z,W)-\frac{1}{n-2}[S(Y,Z)g(X,W)-S(X,Z)g(Y,W)\nonumber\\
&&~~~~~~~~~~~~~~~~~~~~~+g(Y,Z)S(X,W)-g(X,Z)S(Y,W)]\nonumber\\
&&~~~~~~~~~~~~~~~~~~~~~+\frac{r}{(n-1)(n-2)}[g(Y,Z)g(X,W)-g(X,Z)g(Y,W)],
\end{eqnarray}
where $X,Y,Z,W \in \rchi(M),$ $\tilde R$ and $r$ are the Riemannian curvature tensor and scalar curvature of $M$, respectively.\\\\
The $(0,4)$ conharmonic curvature tensor($\cite{jmh61}$) which remains invariant under conharmonic transformation is given by
\begin{eqnarray}\label{2.14}
&&H(X,Y,Z,W)=\tilde R(X,Y,Z,W)-\frac{1}{n-2}[g(X,W)S(Y,Z)-g(X,Z)S(Y,W)\nonumber\\
&&~~~~~~~~~~~~~~~~~~~+g(Y,Z)S(X,W)-g(Y,W)S(X,Z)],
\end{eqnarray}
where $X,Y,Z,W \in \rchi(M)$ and $\tilde R$ is the $(0,4)$ Riemannian curvature tensor.\\\\
The semiconformal curvature tensor has recently been introduced by J. Kim($\cite{jmh64}$, $\cite{jmh65}$) and is given by
\begin{equation}\label{2.15}
\tilde {C}(X,Y,Z,W)=-(n-2)mC(X,Y,Z,W)+[l+(n-2)m]H(X,Y,Z,W),
\end{equation}
where $X,Y,Z,W \in \rchi(M)$, $l,~m$ are nonzero scalar functions. $C$ and $H$ are the conformal and conharmonic curvature tensor respectively. Recently in $2020$ M. Ali, N. A. Pundeer and A. Ali $\cite{jmh66}$ have worked on some applications of semiconformal curvature tensor in general relativity.\\\\
Combining the equations (\ref{2.13}), (\ref{2.14}) and (\ref{2.15}) we get
\begin{eqnarray}\label{2.16}
&&\tilde C(X,Y,Z,W)=l\tilde R(X,Y,Z,W)-\frac{l}{n-2}[S(Y,Z)g(X,W)-S(X,Z)g(Y,W)\nonumber\\
&&~~~~~~~~~~~~~~~~~~~~~+g(Y,Z)S(X,W)-g(X,Z)S(Y,W)]\nonumber\\
&&~~~~~~~~~~~~~~~~~~~~~-\frac{mr}{n-1}[g(Y,Z)g(X,W)-g(X,Z)g(Y,W)].
\end{eqnarray}
$~~~~$Let $R$ be the $(0,3)$ Riemannian curvature tensor of $M$. The $k$-nullity distribution $N(k)$ of a Riemannian manifold $M$ is defined by 
\begin{eqnarray}\label{2.17}
&&N(k):p\rightarrow N_{p}(k)=\nonumber\\
&&\{Z\in T_{p}(M):R(X,Y)Z=k[g(Y,Z)X-g(X,Z)Y]\},
\end{eqnarray}
for all $X,Y \in  T_{p}(M)$, where $k$ is a smooth function. For a quasi-Einstein manifold $M$, if the generator vector field $\xi_{1}$ belongs to some $N(k)$, then $M$ is said to be $N(k)$-quasi Einstein manifold $\cite{jmh67}$.\\\\
$~~~~$The concept of perfect fluid spacetime arose while discussing the structure of this universe. Perfect fluids are often used in the general relativity to model the idealised distribution of matter, such as the interior of a star or isotropic pressure. In general relativity the matter content of the spacetime is described by the energy-momentum tensor. The matter content is assumed to be a fluid having density and pressure and possessing dynamical and kinematical quantities like velocity, acceleration, vorticity, shear and expansion. The energy-momentum tensor $T$ of a perfect fluid spacetime is given by the following equation $\cite{jmh13}$
\begin{equation} \label{2.18}
T(X,Y)= (\sigma+p)A(X)A(Y)+ pg(X,Y).
\end{equation}
Here $g(X,\xi_{1})=A(X), A(\xi_{1})= -1$, for any $X,Y$. $p$ and $\sigma$ are called the isotropic pressure and the energy density respectively. $\xi_{1}$ being the unit timelike velocity vector field.\\\\
The Einstein field equation without cosmological constant $\cite{jmh30}$ is given by
\begin{equation} \label{2.19}
S(X,Y)-\frac{r}{2}g(X,Y)= kT(X,Y);~\forall X,Y \in TM,
\end{equation}
here $r$ being the scalar curvature, $S$ being the Ricci tensor of type $(0,2)$. $k$ is the gravitational constant. From Einstein's field equation it follows that energy momentum tensor is a symmetric $(0,2)$ type tensor of divergence zero.The Einstein equations are fundamental in the construction of cosmological models which imply that the matter determines the geometry of the spacetime and conversely the motion of matter is determined by the metric tensor of the space which is non-flat.
\section{Some geometric properties of $N(QE)_{n}$ and $M(QE)_{n}$ manifolds}\label{s3}
\textbf{Definition 3.1:} {\em Consider a manifold $M^{n}(n \ge 3)$ and a curvature tensor $J$ defined on it. Then the manifold is called $J$-flat if $J(X,Y)Z=0$, for all $X, Y, Z \in \rchi(M).$}\\\\
$~~~~$Consider the $N(QE)_{n}$ manifold defined by (\ref{2.3}) and let it is semi-conformally flat. Then by the definition $(3.1)$ we get, $\tilde C(X,Y)Z=0$, which implies $\tilde C(X,Y,Z,W)=0.$ Thus from (\ref{2.16}) for all $X, Y, Z, W \in \rchi(M)$ we get
\begin{eqnarray}\label{3.1}
&&\tilde R(X,Y,Z,W)=\frac{1}{n-2}[S(Y,Z)g(X,W)-S(X,Z)g(Y,W)\nonumber\\
&&~~~~~~~~~~~~~~~~~~~~~+g(Y,Z)S(X,W)-g(X,Z)S(Y,W)]\nonumber\\
&&~~~~~~~~~~~~~~~~~~~~~+\frac{mr}{l(n-1)}[g(Y,Z)g(X,W)-g(X,Z)g(Y,W)].
\end{eqnarray}
Using the equation (\ref{2.3}) we have
\begin{eqnarray}\label{3.2}
&&\tilde R(X,Y,Z,W)=\{\frac{mr}{l(n-1)}+\frac{2\alpha}{n-2}\}[g(Y,Z)g(X,W)-g(X,Z)g(Y,W)]\nonumber\\
&&~~~~~~~~~~~~~~~~~~~~~+\frac{\beta}{n-2}[D(Y,Z)g(X,W)-D(X,Z)g(Y,W)\nonumber\\
&&~~~~~~~~~~~~~~~~~~~~~+g(Y,Z)D(X,W)-g(X,Z)D(Y,W)].
\end{eqnarray}

From the equation (\ref{2.8}) we realise that this is a manifold of nearly quasi-constant curvature. Hence we derive the following theorem as:\\\\
\textbf{Theorem 3.1:} {\em A semiconformally flat $N(QE)_{n}$ manifold is an $N(QC)_{n}$.}\\\\
Using the remark (2.1) we state the following corollary:\\\\
\textbf{Corollary 3.1:} {\em A semiconformally flat $M(QE)_{n}$ manifold is an $N(QC)_{n}$.}\\\\
Now from the equations (\ref{2.5}) and (\ref{3.1}) we get
\begin{eqnarray}\label{3.3}
&&\tilde R(X,Y,Z,W)=\{\frac{mr}{l(n-1)}+\frac{2\alpha}{n-2}\}[g(Y,Z)g(X,W)-g(X,Z)g(Y,W)]\nonumber\\
&&~~~~~~~~~~~~~~~~~~~~~+\frac{\beta}{n-2}[g(X,W)\{A(Y)B(Z)+A(Z)B(Y)\}\nonumber\\
&&~~~~~~~~~~~~~~~~~~~~~-g(Y,W)\{A(X)B(Z)+A(Z)B(X)\}\nonumber\\
&&~~~~~~~~~~~~~~~~~~~~~+g(Y,Z)\{A(X)B(W)+A(W)B(X)\}\nonumber\\
&&~~~~~~~~~~~~~~~~~~~~~-g(X,Z)\{A(Y)B(W)+A(W)B(Y)\}].
\end{eqnarray}
And hence from (\ref{2.9}) we realise that this is a manifold of mixed quasi-constant curvature. Hence we derive the following theorem as\\\\
\textbf{Theorem 3.2:} {\em A semiconformally flat $M(QE)_{n}$ manifold is an $M(QC)_{n}.$}\\\\
Now consider the $N(QE)_{n}$ manifold is conharmonically flat. Then, $H(X,Y,Z,W)=0$, hence from the equation (\ref{2.14}) we get
\begin{eqnarray}\label{3.4}
&&\tilde R(X,Y,Z,W)=\frac{1}{n-2}[g(X,W)S(Y,Z)-g(X,Z)S(Y,W)\nonumber\\
&&~~~~~~~~~~~~~~~~~~~+g(Y,Z)S(X,W)-g(Y,W)S(X,Z)].
\end{eqnarray}
Using the equation (\ref{2.3}) we get
\begin{eqnarray}\label{3.5}
&&\tilde R(X,Y,Z,W)=\frac{2\alpha}{n-2}[g(Y,Z)g(X,W)-g(X,Z)g(Y,W)]\nonumber\\
&&~~~~~~~~~~~~~~~~~~~~~~+\frac{\beta}{n-2}[D(Y,Z)g(X,W)-D(X,Z)g(Y,W)\nonumber\\
&&~~~~~~~~~~~~~~~~~~~~~+g(Y,Z)D(X,W)-g(X,Z)D(Y,W)].
\end{eqnarray}
From the equation (\ref{2.8}) we realise that this is a manifold of nearly quasi-constant curvature. Hence we derive the following theorem as:\\\\
\textbf{Theorem 3.3:} {\em A conharmonically flat $N(QE)_{n}$ manifold is an $N(QC)_{n}$.}\\\\
Using the remark (2.1) we state the following corollary:\\\\
\textbf{Corollary 3.2:} {\em A conharmonically flat $M(QE)_{n}$ manifold is an $N(QC)_{n}$.}\\\\
Again considering (\ref{2.5}) and (\ref{3.4}) we get
\begin{eqnarray}\label{3.6}
&&\tilde R(X,Y,Z,W)=\frac{2\alpha}{n-2}[g(Y,Z)g(X,W)-g(X,Z)g(Y,W)]\nonumber\\
&&~~~~~~~~~~~~~~~~~~~~~~\frac{\beta}{n-2}[g(X,W)\{A(Y)B(Z)+A(Z)B(Y)\}\nonumber\\
&&~~~~~~~~~~~~~~~~~~~~~-g(Y,W)\{A(X)B(Z)+A(Z)B(X)\}\nonumber\\
&&~~~~~~~~~~~~~~~~~~~~~+g(Y,Z)\{A(X)B(W)+A(W)B(X)\}\nonumber\\
&&~~~~~~~~~~~~~~~~~~~~~-g(X,Z)\{A(Y)B(W)+A(W)B(Y)\}].
\end{eqnarray}
Thus from (\ref{2.9}) we can say that this manifold is an $M(QC)_{n}$. So, the following theorem is stated:\\\\
\textbf{Theorem 3.4:} {\em A conharmonically flat $M(QE)_{n}$ manifold is an $M(QC)_{n}$.}\\\\
Now consider the $N(QE)_{n}$ manifold is $W_{2}$-flat. Then, $W_{2}(X,Y,Z,W)=0$, hence from the equation (\ref{2.11}) we get
\begin{equation} \label{3.7}
\tilde R(X,Y,Z,W)=\frac{1}{n-1}[g(Y,Z)S(X,W)-g(X,Z)S(Y,W)].
\end{equation}
Using the equation (\ref{2.3}) we get
\begin{eqnarray} \label{3.8}
&&\tilde R(X,Y,Z,W)=\frac{\alpha }{n-1}[g(Y,Z)g(X,W)-g(X,Z)g(Y,W)]\nonumber\\
&&~~~~~~~~~~~~~~~~~~~+\frac{\beta}{n-1}[g(Y,Z)D(X,W)-g(X,Z)D(Y,W)].
\end{eqnarray}
From the equation (\ref{2.8}) we realise that this is not a manifold of nearly quasi-constant curvature. Hence we derive the following theorem as:\\\\
\textbf{Theorem 3.5:} {\em A $W_{2}$-flat $N(QE)_{n}$ manifold is NOT an $N(QC)_{n}$.}\\\\
From the remark (2.1) we see that a an $M(QE)_{n}$ manifold is an $N(QE)_{n}$ manifold. Again from the remark (2.2) we see that a manifold which is not an $N(QC)_{n}$ cannot be an $M(QC)_{n}$. So, the following corollary is stated:\\\\
\textbf{Corollary 3.3} {\em A $W_{2}$-flat $M(QE)_{n}$ manifold is NOT an $M(QC)_{n}$.}\\\\
Now consider the $N(QE)_{n}$ manifold is projectively flat. Then, $P(X,Y,Z,W)=0$, hence from the equation (\ref{2.12}) we get
\begin{equation}\label{3.7}
\tilde R(X,Y,Z,W)=\frac{1}{n-1}\{S(Y,Z)g(X,W)-S(X,Z)g(Y,W)\}.
\end{equation}
From the equation (\ref{2.3}) we get
\begin{eqnarray} \label{3.8}
&&\tilde R(X,Y,Z,W)=\frac{\alpha }{n-1}[g(Y,Z)g(X,W)-g(X,Z)g(Y,W)]\nonumber\\
&&~~~~~~~~~~~~~~~~~~~+\frac{\beta}{n-1}[D(Y,Z)g(X,W)-D(X,Z)g(Y,W)].
\end{eqnarray}
From the equation (\ref{2.8}) we realise that this is not a manifold of nearly quasi-constant curvature. Hence we derive the following theorem as:\\\\
\textbf{Theorem 3.6:} {\em A projectively flat $N(QE)_{n}$ manifold is NOT an $N(QC)_{n}$.}\\\\
\textbf{Corollary 3.4:} {\em A projectively flat $M(QE)_{n}$ manifold is NOT an $M(QC)_{n}$.}\\\\
\textbf{Remark 3.1:} {\em A $J$-flat $N(QE)_{n}$ manifold is not necessarily an $N(QC)_{n}$, and a $J$-flat $M(QE)_{n}$ manifold is not necessarily an $M(QC)_{n}$ where $J$ is some curvature tensor defined on the manifold.}
\section{$W_{2}$-Ricci pseudosymmetric $M(QE)_{4}$ spacetime} \label{s4}
\textbf{Definition 4.1:} $\cite{jmh59}$ {\em A semi-Riemannian manifold $(M^{n},g), (n \ge 3)$ admitting a $(0,4)$ tensor field $J$ is said to be $J$-Ricci pseudosymmetric if $J\cdot S$ and $Q(g, S)$ are linearly dependent, i.e., $J \cdot S=L_{S}Q(g, S)$ holds on a set $U_{S}=\{x\in M :S \neq \frac{r}{n}g$ at $x\}$, where $L_{S}$ is some function on $U_{S}.$}\\\\
If $J=W_{2}$ or $P$ then the manifold is called $W_{2}$-Ricci pseudosymmetric or projectively Ricci pseudosymmetric respectively.\\\\
$~~~~$In this section a $W_{2}$-Ricci pseudosymmetric $M(QE)_{4}$ spacetime is considered, i.e., we take an $M(QE)_{4}$ spacetime satisfying the condition $W_{2}. S = F_{S} Q(g,S)$. Here $F_{S}$ being a certain function on the set $U_{S} = \{x\in M :S \neq \frac{r}{n}g$ at $x\}$ and $Q(g,S)$ being the Tachibana tensor working on the metric tensor and the Ricci tensor. Now for all $ X,Y,Z,W \in \rchi(M^4)$, we have
\begin{eqnarray} \label{4.1}
& & S(W_{2}(X,Y)Z,W)+ S(Z,W_{2}(X,Y)W)\nonumber\\
&&= F_{S}[g(Y,Z)S(X,W)-g(X,Z)S(Y,W)\nonumber\\
&&+ g(Y,W)S(Z,X)- g(X,W)S(Y,Z)].
\end{eqnarray}
Using the equation (\ref{2.5}) from the equation (\ref{4.1}) we get
\begin{eqnarray} \label{4.2}
&&\alpha g(W_{2}(X,Y)Z,W)+\beta[A(W_{2}(X,Y)Z)B(W)\nonumber\\
&&+B(W_{2}(X,Y)Z)A(W)]+\alpha g(W_{2}(X,Y)W,Z)\nonumber\\
&&+\beta [A(W_{2}(X,Y)W)B(Z)+B(W_{2}(X,Y)W)A(Z)] \nonumber\\
&&=F_{S}[g(Y,Z)S(X,W)-g(X,Z)S(Y,W) \nonumber\\
&&~~+g(Y,W)S(Z,X)-g(X,W)S(Y,Z)].
\end{eqnarray}
Using $g(W_{2}(X,Y)Z,W)=W_{2}(X,Y,Z,W)$, equation (\ref{4.2}) can be written as
\begin{eqnarray} \label{4.3}
&&\alpha W_{2}(X,Y,Z,W)+\beta[W_{2}(X,Y,Z,\xi_{1})B(W)\nonumber\\
&&+W_{2}(X,Y,Z,\xi_{2})A(W)]+\alpha W_{2}(X,Y,W,Z)\nonumber\\
&&+\beta [W_{2}(X,Y,W,\xi_{1})B(Z)+W_{2}(X,Y,W,\xi_{2})A(Z)] \nonumber\\
&&=F_{S}[g(Y,Z)S(X,W)-g(X,Z)S(Y,W) \nonumber\\
&&~~+g(Y,W)S(Z,X)-g(X,W)S(Y,Z)].
\end{eqnarray}
Putting $Z=W=\xi_{1},~n=4$ and using the equations (\ref{2.5}), (\ref{2.11}) from the equation (\ref{4.3}) we get
\begin{eqnarray}\label{4.4}
&&\frac{-\alpha \beta}{3}[A(X)B(Y)-A(Y)B(X)]\nonumber\\
&&-\beta[\tilde R(X,Y,\xi_{1},\xi_{2})+\frac{\alpha}{3}\{A(X)B(Y)-A(Y)B(X)\}]\nonumber\\
&&=\beta F_{S}[A(X)B(Y)-A(Y)B(X)].
\end{eqnarray}
Since $\beta \neq 0$ thus from (\ref{4.4}) we get
\begin{equation}\label{4.5}
\tilde R(X,Y,\xi_{1},\xi_{2})=(-\frac{2\alpha}{3}-F_{S})[A(X)B(Y)-A(Y)B(X)].
\end{equation}
Again after contraction over $X, W$ and  putting $Z=\xi_{1}$ from the equation (\ref{4.3}) we get
\begin{equation}\label{4.6}
-\frac{5\alpha \beta}{3}+\frac{\beta^{2}}{3}A(Y)+\beta \tilde R(\xi_{1},Y,\xi_{1},\xi_{2})=4\beta B(Y)F_{S}.
\end{equation}
Since $\beta \neq 0$ thus by using the equation (\ref{4.5}) from (\ref{4.6}) we obtain
\begin{equation}\label{4.7}
\frac{\beta}{3}A(Y)+(-\alpha-3F_{S})B(Y)=0.
\end{equation}
Putting $Y=\xi_{2}$ from (\ref{4.7}) we get
\begin{equation}\label{4.8}
F_{S}=-\frac{\alpha}{3}.
\end{equation}
From (\ref{4.5}) and (\ref{4.8}) we get
\begin{equation}\label{4.9}
\tilde R(X,Y,\xi_{1},\xi_{2})=\frac{\alpha}{3}[A(Y)B(X)-A(X)B(Y)].
\end{equation}
Equation (\ref{4.9}) can be written as
\begin{equation}\label{4.10}
R(X,Y)\xi_{1}=\frac{\alpha}{3}[A(Y)X-A(X)Y].
\end{equation}
Thus from the equations (\ref{2.17}) and (\ref{4.10}) we realise that for this manifold the generator vector field belongs to the $k$-nullity distribution for $k=\frac{\alpha}{3}$ and hence is an $N(\frac{\alpha}{3})$ quasi-Einstein spacetime. So, we can derive the following theorem as:\\\\
\textbf{Theorem 4.1:} {\em A $W_{2}$-Ricci pseudosymmmetric $M(QE)_{4}$ spacetime is an $N(\frac{\alpha}{3})$-quasi Einstein spacetime.}\\\\
Taking contraction with respect to $X$ from the equation (\ref{4.10}) we get
\begin{equation}\label{4.11}
S(Y,\xi_{1})=\alpha A(Y).
\end{equation}
Equation (\ref{4.11}) can be written as
\begin{equation}\label{4.12}
QY=\alpha Y.
\end{equation}
Hence we derive the following theorem:\\\\
\textbf{Theorem 4.2:} {\em In a $W_{2}$-Ricci pseudosymmmetric $M(QE)_{4}$ spacetime the eigenvalue of the Ricci operator $Q$ is $\alpha$.}
\section{Projectively Ricci pseudosymmetric $M(QE)_{4}$ spacetime} \label{s5}
$~~~~$In this section a projectively Ricci pseudosymmetric $M(QE)_{4}$ spacetime is considered, i.e., we take an $M(QE)_{4}$ spacetime satisfying the condition $P. S = F_{S} Q(g,S)$. Here $F_{S}$ being a certain function on the set $U_{S} = \{x\in M :S\neq \frac{r}{n}g$ at $x\}$ and $Q(g,S)$ being the Tachibana tensor working on the metric tensor and the Ricci tensor. Now for all $ X,Y,Z,W \in \rchi(M^4)$, we have
\begin{eqnarray} \label{5.1}
& & S(P(X,Y)Z,W)+ S(Z,P(X,Y)W)\nonumber\\
&&= F_{S}[g(Y,Z)S(X,W)-g(X,Z)S(Y,W)\nonumber\\
&&+ g(Y,W)S(Z,X)- g(X,W)S(Y,Z)].
\end{eqnarray}
Using the equation (\ref{2.5}) from the equation (\ref{5.1}) we get
\begin{eqnarray} \label{5.2}
&&\alpha g(P(X,Y)Z,W)+\beta[A(P(X,Y)Z)B(W)\nonumber\\
&&+B(P(X,Y)Z)A(W)]+\alpha g(P(X,Y)W,Z)\nonumber\\
&&+\beta [A(P(X,Y)W)B(Z)+B(P(X,Y)W)A(Z)] \nonumber\\
&&=F_{S}[g(Y,Z)S(X,W)-g(X,Z)S(Y,W) \nonumber\\
&&~~+g(Y,W)S(Z,X)-g(X,W)S(Y,Z)].
\end{eqnarray}
Using $g(P(X,Y)Z,W)=P(X,Y,Z,W)$, equation (\ref{5.2}) can be written as
\begin{eqnarray} \label{5.3}
&&\alpha P(X,Y,Z,W)+\beta[P(X,Y,Z,\xi_{1})B(W)\nonumber\\
&&+P(X,Y,Z,\xi_{2})A(W)]+\alpha P(X,Y,W,Z)\nonumber\\
&&+\beta [P(X,Y,W,\xi_{1})B(Z)+P(X,Y,W,\xi_{2})A(Z)] \nonumber\\
&&=F_{S}[g(Y,Z)S(X,W)-g(X,Z)S(Y,W) \nonumber\\
&&~~+g(Y,W)S(Z,X)-g(X,W)S(Y,Z)].
\end{eqnarray}
Putting $Z=W=\xi_{1}$ we get
\begin{equation}\label{5.4}
\alpha P(X,Y,\xi_{1},\xi_{1})-\beta P(X,Y,\xi_{1},\xi_{2})=-\beta F_{S}[A(Y)B(X)-A(X)B(Y)].
\end{equation}
Using (\ref{2.12}), from (\ref{5.4}) we get
\begin{equation}\label{5.5}
-\beta \tilde R(X,Y,\xi_{1},\xi_{2})=-\beta F_{S}[A(Y)B(X)-A(X)B(Y)].
\end{equation}
$\beta \neq 0$ implies
\begin{equation}\label{5.6}
\tilde R(X,Y,\xi_{1},\xi_{2})=F_{S}[A(Y)B(X)-A(X)B(Y)].
\end{equation}
Now taking the contraction over $X$ and $W$ and using the equation (\ref{2.6}) from the equation (\ref{5.3}) we get
\begin{eqnarray}\label{5.7}
\beta ^{2}A(Y)+\beta \tilde R(\xi_{1},Y,\xi_{1},\xi_{2})+\alpha \beta B(Y)=4\beta F_{S}B(Y).
\end{eqnarray}
Using the equation (\ref{5.6}) we get
\begin{equation}\label{5.8}
\beta ^{2}A(Y)+\alpha \beta B(Y)=3\beta F_{S}B(Y).
\end{equation}
Since $\beta \neq 0$ thus putting $Y=\xi_{2}$ in (\ref{5.8}) we get
\begin{equation}\label{5.9}
F_{S}=\frac{\alpha}{3}.
\end{equation}
Putting the value of $F_{S}$ in the equation (\ref{5.6}) we get
\begin{equation}\label{5.10}
\tilde R(X,Y,\xi_{1},\xi_{2})=\frac{\alpha}{3}[A(Y)B(X)-A(X)B(Y)].
\end{equation}
Which can be written as
\begin{equation}\label{5.11}
R(X,Y)\xi_{1}=\frac{\alpha}{3}[A(Y)X-A(X)Y].
\end{equation}
Thus from the equations (\ref{2.17}) and (\ref{5.11}) we realise that for this manifold the generator vector field belongs to the $k$-nullity distribution for $k=\frac{\alpha}{3}$ and hence is an $N(\frac{\alpha}{3})$ quasi-Einstein spacetime. So, we can derive the following theorem as:\\\\
\textbf{Theorem 5.1:} {\em A projectively Ricci pseudosymmmetric $M(QE)_{4}$ spacetime is an $N(\frac{\alpha}{3})$-quasi Einstein spacetime.}\\\\
Taking contraction with respect to $X$ from the equation (\ref{5.11}) we get
\begin{equation}\label{5.12}
S(Y,\xi_{1})=\alpha A(Y).
\end{equation}
Equation (\ref{5.12}) can be written as
\begin{equation}\label{5.13}
QY=\alpha Y.
\end{equation}
Hence we derive the following theorem:\\\\
\textbf{Theorem 5.2:} {\em In a projectively Ricci pseudosymmmetric $M(QE)_{4}$ spacetime the eigenvalue of the Ricci operator $Q$ is $\alpha$.}
\section{$M(QE)_{4}$ spacetime satisfying some conditions} \label{s6}
$~~~~$Here we consider an $M(QE)_{4}$ spacetime satisfying the conditions like $H \cdot S=0$ and $\tilde C \cdot S=0$. Firstly, let us consider that the $M(QE)_{4}$ spacetime satisfies the condition $H \cdot S=0$. This gives us the following equation
\begin{eqnarray} \label{6.1}
S(H(X,Y)Z,W)+ S(Z,H(X,Y)W)=0,
\end{eqnarray}
 for all $ X,Y,Z,W \in \rchi(M^4)$. Using the equation (\ref{2.5}) from the equation (\ref{6.1}) we get
\begin{eqnarray} \label{6.2}
&&\alpha g(H(X,Y)Z,W)+\beta[A(H(X,Y)Z)B(W)\nonumber\\
&&+B(H(X,Y)Z)A(W)]+\alpha g(H(X,Y)W,Z)\nonumber\\
&&+\beta [A(H(X,Y)W)B(Z)+B(H(X,Y)W)A(Z)]=0.
\end{eqnarray}
Using $g(H(X,Y)Z,W)=H(X,Y,Z,W)$, equation (\ref{6.2}) can be written as
\begin{eqnarray} \label{6.3}
&&\alpha H(X,Y,Z,W)+\beta[H(X,Y,Z,\xi_{1})B(W)\nonumber\\
&&+H(X,Y,Z,\xi_{2})A(W)]+\alpha H(X,Y,W,Z)\nonumber\\
&&+\beta [H(X,Y,W,\xi_{1})B(Z)+H(X,Y,W,\xi_{2})A(Z)]=0.
\end{eqnarray}
Putting $Z=W=\xi_{1}$ in (\ref{6.3}) and using the equation (\ref{2.14}) we get
\begin{equation}\label{6.4}
\beta \tilde R(X,Y,\xi_{1},\xi_{2})=\beta \alpha [A(Y)B(X)-A(X)B(Y)].
\end{equation}
$\beta \neq 0$ will imply
\begin{equation}\label{6.5}
\tilde R(X,Y,\xi_{1},\xi_{2})=\alpha [A(Y)B(X)-A(X)B(Y)].
\end{equation}
Again after taking contraction over $X,~W$ and putting $Z=\xi_{1}$ from (\ref{6.3}) we get
\begin{equation}\label{6.6}
\beta \tilde R(\xi_{1},Y,\xi_{1},\xi_{2})-3\alpha \beta B(Y)=0.
\end{equation}
Using the equation (\ref{6.5}), equation (\ref{6.6}) gives
\begin{equation}\label{6.7}
\alpha \beta B(Y)=0.
\end{equation}
Putting $Y=\xi_{2}$, we get
\begin{equation}\label{6.8}
\alpha \beta=0.
\end{equation}
Since $\beta \neq 0$ thus from (\ref{6.8}) we get $\alpha=0$. From the equation (\ref{2.6}) we get $r=0,$ i.e. the scalar curvature of the manifold is $0$. Thus we can state the following theorem as:\\\\
\textbf{Theorem 6.1:} {\em The scalar curvature of an $M(QE)_{4}$ spacetime satisfying the condition $H \cdot S=0~is~0$.}\\\\
Next we consider an $M(QE)_{4}$ spacetime satisfying the condition $\tilde C \cdot S=0.$ This gives us the equation
\begin{eqnarray} \label{6.9}
S(\tilde C(X,Y)Z,W)+ S(Z,\tilde C(X,Y)W)=0,
\end{eqnarray}
 for all $ X,Y,Z,W \in \rchi(M^4)$. Using the equation (\ref{2.5}) from the equation (\ref{6.9}) we get
\begin{eqnarray} \label{6.10}
&&\alpha g(\tilde C(X,Y)Z,W)+\beta[A(\tilde C(X,Y)Z)B(W)\nonumber\\
&&+B(\tilde C(X,Y)Z)A(W)]+\alpha g(\tilde C(X,Y)W,Z)\nonumber\\
&&+\beta [A(\tilde C(X,Y)W)B(Z)+B(\tilde C(X,Y)W)A(Z)]=0.
\end{eqnarray}
Using $g(\tilde C(X,Y)Z,W)=\tilde C(X,Y,Z,W)$, equation (\ref{6.10}) can be written as
\begin{eqnarray} \label{6.11}
&&\alpha \tilde C(X,Y,Z,W)+\beta[\tilde C(X,Y,Z,\xi_{1})B(W)\nonumber\\
&&+\tilde C(X,Y,Z,\xi_{2})A(W)]+\alpha \tilde C(X,Y,W,Z)\nonumber\\
&&+\beta [\tilde C(X,Y,W,\xi_{1})B(Z)+\tilde C(X,Y,W,\xi_{2})A(Z)]=0.
\end{eqnarray}
Putting $Z=W=\xi_{1}$ in (\ref{6.11}) and using (\ref{2.16}) we get
\begin{equation}\label{6.12}
l\beta \tilde R(X,Y,\xi_{1},\xi_{2})=l\beta \alpha [A(Y)B(X)-A(X)B(Y)].
\end{equation}
$l \neq 0$ and $\beta \neq 0$ will imply
\begin{equation}\label{6.13}
\tilde R(X,Y,\xi_{1},\xi_{2})=\alpha [A(Y)B(X)-A(X)B(Y)].
\end{equation}
Again after contracting (\ref{6.11}) over $X$ and $W$ and putting $Z=\xi_{1}$ we get
\begin{eqnarray}\label{6.14}
l[\beta \tilde R(\xi_{1},Y,\xi_{1},\xi_{2})-3\alpha \beta B(Y)]-4m\alpha A(Y)=0.
\end{eqnarray}
Using the equation (\ref{6.13}) we get
\begin{equation}\label{6.15}
-2l\alpha \beta B(Y)-4m\alpha A(Y)=0.
\end{equation}
Putting $Y=\xi_{1}$ we get from (\ref{6.15}),
\begin{equation}\label{6.16}
m\alpha=0.
\end{equation}
Since $m \neq 0$ thus from (\ref{6.16}) we get $\alpha=0$. From the equation (\ref{2.6}) we get $r=0,$ i.e. the scalar curvature of the manifold is $0$. Again putting $Y=\xi_{2}$ from the equation (\ref{6.15}) we get
\begin{equation}\label{6.17}
l\alpha \beta=0.
\end{equation}
Since $\beta \neq 0$ and $l\neq 0$ thus from (\ref{6.17}) we get $\alpha=0$. From the equation (\ref{2.6}) we get $r=0,$ i.e. the scalar curvature of the manifold is 0. So, in both the cases we see that the scalar curvature of the manifold is 0. Hence we can derive the following theorem as:\\\\
\textbf{Theorem 6.2:} {\em The scalar curvature of an $M(QE)_{4}$ spacetime satisfying the condition $\tilde C \cdot S=0$ is $0$.}
\section{Semiconformally flat Ricci pseudosymmetric\\$M(QE)_{4}$ spacetime} \label{s7}
$~~~~$In this section we study the geometric nature of a Ricci pseudosymmetric $M(QE)_{4}$ spacetime. Consider an $M(QE)_{4}$ spacetime and suppose it is Ricci pseudosymmetric. Then we have the following condition satisfied,
\begin{equation}\label{7.1}
(R(X,Y) \cdot S(Z,W)) = F_{S} Q(g,S)(Z,W;X,Y).
\end{equation}
Here $X,Y,Z,W \in \rchi(M),$ $F_{S}$ being a certain function on the set $U_{S} = \{x\in M : S \neq \frac{r}{n}g$ at $x\}$ and $Q(g,S)$ being the Tachibana tensor working on the metric tensor and the Ricci tensor. The LHS of the equation (\ref{7.1}) can be expressed as
\begin{equation}\label{7.2}
(R(X,Y) \cdot S(Z,W)) = -S(R(X,Y)Z,W)-S(Z,R(X,Y)W)
\end{equation}
and the RHS of the equation (\ref{7.1}) can be expressed as
\begin{equation}\label{7.3}
F_{S} Q(g,S)(Z,W;X,Y)=-F_{S}[S((X\Lambda_{g}Y)Z,W)+S(Z,(X\Lambda_{g}Y)W)],
\end{equation}
where,
\begin{equation}\label{7.4}
(X\Lambda_{g}Y)Z=g(Y,Z)X-g(X,Z)Y.
\end{equation}
Using (\ref{7.2}), (\ref{7.3}) and (\ref{7.4}) from the equation (\ref{7.1}) we get
\begin{eqnarray}\label{7.5}
&&S(R(X,Y)Z,W)+S(Z,R(X,Y)W)=\nonumber\\
&&F_{S}[g(Y,Z)S(X,W)-g(X,Z)S(Y,W)\nonumber\\
&&+g(Y,W)S(X,Z)-g(X,W)S(Y,Z)].
\end{eqnarray}
With the help of the equation (\ref{2.5}) the equation (\ref{7.5}) takes the form
\begin{eqnarray}\label{7.6}
&&\alpha g(R(X,Y)Z,W)+\alpha g(Z,R(X,Y)W)\nonumber\\
&&\beta[A(R(X,Y)Z)B(W)+A(W)B(R(X,Y)Z)\nonumber\\
&&+A(Z)B(R(X,Y)W)+A(R(X,Y)W)B(Z)]\nonumber\\
&&=F_{S}[g(Y,Z)\{\alpha g(X,W)+\beta [A(X)B(W)+A(W)B(X)]\}\nonumber\\
&&-g(X,Z)\{\alpha g(Y,W)+\beta [A(Y)B(W)+A(W)B(Y)]\}\nonumber\\
&&+g(Y,W)\{\alpha g(X,Z)+\beta [A(X)B(Z)+A(Z)B(X)]\}\nonumber\\
&&-g(X,W)\{\alpha g(Y,Z)+\beta [A(Y)B(Z)+A(Z)B(Y)]\}].
\end{eqnarray}
Putting $Z=W=\xi_{1}$ in (\ref{7.6}) and using the condition $\beta \neq 0$ we get
\begin{equation}\label{7.7}
\tilde R(X,Y,\xi_{1},\xi_{2})=F_{S}[A(Y)B(X)-A(X)B(Y)].
\end{equation}
Again putting $Z=\xi_{1}, W=\xi_{2}$ in the equation (\ref{3.3}) and using $n=4, r=4\alpha$ we get
\begin{equation}\label{7.8}
\tilde R(X,Y,\xi_{1},\xi_{2})=\alpha(1+\frac{4m}{3l})[A(Y)B(X)-A(X)B(Y)].
\end{equation}
From (\ref{7.7}) and (\ref{7.8}) we get
\begin{equation}\label{7.9}
\{F_{S}-\alpha(1+\frac{4m}{3l})\}[A(Y)B(X)-A(X)B(Y)]=0.
\end{equation}
So, either $A(Y)B(X)-A(X)B(Y)=0$ or $F_{S}-\alpha(1+\frac{4m}{3l})=0$. Considering $A(Y)B(X)-A(X)B(Y)\neq 0$ we must get, $F_{S}-\alpha(1+\frac{4m}{3l})=0$. Hence from the equation (\ref{7.7}) we get
\begin{equation}\label{7.10}
\tilde R(X,Y,\xi_{1},\xi_{2})=\alpha(1+\frac{4m}{3l})[A(Y)B(X)-A(X)B(Y)],
\end{equation}
which can be written as
\begin{equation}\label{7.11}
R(X,Y)\xi_{1}=\alpha(1+\frac{4m}{3l})[A(Y)X-A(X)Y].
\end{equation}
Using the equation (\ref{2.17}) from the equation (\ref{7.11}) we realise that the vector field $\xi_{1}$ belongs to $k$-nullity distribution for $k=\alpha(1+\frac{4m}{3l})$. Hence we state the following theorem:\\\\
\textbf{Theorem 7.1:} {\em A semiconformally flat Ricci pseudosymmetric $M(QE)_{4}$ spacetime with the condition $A(Y)B(X)-A(X)B(Y)\neq 0$, is an $N(\alpha(1+\frac{4m}{3l}))$-quasi Einstein spacetime.}
\section{Perfect fluid $M(QE)_{4}$ spacetime satisfying Einstein field equation} \label{s8}
$~~~~$In this section we discuss about a perfect fluid $M(QE)_{4}$ spacetime satisfying Einstein field equation without cosmological constant. Consider a $M(QE)_{4}$ perfect fluid spacetime satisfying Einstein field equation without cosmological constant. Then form the equations (\ref{2.18}) and (\ref{2.19}) we get
\begin{equation}\label{8.1}
S(X,Y)=k(\sigma +p)A(X)A(Y)+(pk+\frac{r}{2})g(X,Y).
\end{equation}
Putting $X=Y=\xi_{1}$ in (\ref{8.1}) and using the equation (\ref{2.5}) we get
\begin{equation}\label{8.2}
-\alpha=k\sigma -\frac{r}{2}.
\end{equation}
Using the equation (\ref{2.6}), from (\ref{8.2}) we get
\begin{equation}\label{8.3}
\alpha=k \sigma.
\end{equation}
From (\ref{2.6}) and (\ref{8.3}) we get
\begin{equation}\label{8.4}
r=4k \sigma.
\end{equation}
Now, replacing $X$ by $QX$ from the equation (\ref{8.1}) we have
\begin{equation}\label{8.5}
S(QX,Y)=k(\sigma+p)S(X,\xi_{1})g(Y,\xi_{1})+(pk+\frac{r}{2})S(X,Y).
\end{equation}
Taking contraction over $X,Y$ from the equation (\ref{8.5}) we get
\begin{equation}\label{8.6}
||Q||^{2}=-k\alpha(\sigma+p)+(pk+\frac{r}{2})r.
\end{equation}
Using (\ref{8.3}) and (\ref{8.4}) from the equation (\ref{8.6}) we get
\begin{equation}\label{8.7}
||Q||^{2}=k^{2}\sigma(3p+7\sigma).
\end{equation}
Now for a perfect fluid $M(QE)_{4}$ spacetime satisfying Einstein field equation without cosmological constant  it is proved($\cite{jmh69}$) that
\begin{equation}\label{8.8}
\sigma+p=0.
\end{equation}
They showed that such a spacetime represents inflation $\cite{jmh73}$ and also termed as Phantom Barrier $\cite{jmh74}$. If the perfect fluid spacetime satisfies the condition (\ref{8.8}) then the fluid behaves like a cosmological constant $\cite{jmh72}$. Now combining the equations (\ref{8.7}), (\ref{8.8}) and using (\ref{8.4}) we get
\begin{equation}\label{8.9}
||Q||^{2}=\frac{r^{2}}{4}.
\end{equation}
Hence we can derive the following theorem:\\\\
\textbf{Theorem 8.1} {\em In a perfect fluid $M(QE)_{4}$ spacetime satisfying Einstein field equation without cosmological constant the square of the length of the Ricci operator $Q$ is $\frac{r^{2}}{4}$, where $r$ is the scalar curvature of the spacetime.}\\\\
Now if the spacetime is flat then $r=0$. Then from the equation (\ref{8.4}) we get
\begin{equation}\label{8.10}
k\sigma=0.
\end{equation}
Since $k \neq 0,$ thus we must have $\sigma=0$. Hence from (\ref{8.8}) we get, $p=0$. Thus the spacetime becomes a vacuum. Hence we can derive the following theorem:\\\\
\textbf{Theorem 8.2} {\em A flat perfect fluid $M(QE)_{4}$ spacetime satisfying Einstein field equation without cosmological constant is a vacuum.}\\\\
Now, let us consider that $\xi_{1}$ is a Killing vector field. Then, we must get the equation
\begin{equation}\label{8.11}
(\mathcal{L}_{\xi_{1}}g)(X,Y)=0,
\end{equation}
for all $X, Y \in \rchi(M)$. Again, using the equations (\ref{2.5}) and (\ref{2.19}) we get
\begin{equation}\label{8.12}
(\alpha-\frac{r}{2}) g(X,Y)+\beta [A(X)B(Y)+A(Y)B(X)]=kT(X,Y).
\end{equation}
If $\alpha, \beta$ are constants then from (\ref{2.6}) $r$ is also constant. Thus taking the Lie derivative to both the sides of the equation (\ref{8.12}) we get
\begin{equation}\label{8.13}
(\alpha-\frac{r}{2})(\mathcal{L}_{\xi_{1}}g)(X,Y)=k(\mathcal{L}_{\xi_{1}}T)(X,Y).
\end{equation}
Using (\ref{8.11}) we get
\begin{equation}\label{8.14}
k(\mathcal{L}_{\xi_{1}}T)(X,Y)=0.
\end{equation}
$k \neq 0$ yields
\begin{equation}\label{8.15}
(\mathcal{L}_{\xi_{1}}T)(X,Y)=0.
\end{equation}
Thus we see that the Lie derivative of the energy-momentum tensor vanishes. Conversely let the Lie derivative of the energy-momentum tensor vanishes. So if $\alpha, \beta$ are constants then taking Lie derivative to both the sides of the equation (\ref{8.12}) we observe that $(\alpha-\frac{r}{2})(\mathcal{L}_{\xi_{1}}g)(X,Y)=0.$ Now if we consider $\alpha \neq 0$ then using (\ref{2.6}) we get $\alpha-\frac{r}{2}=-\alpha \neq 0.$ Hence we derive the following theorem:\\\\
\textbf{Theorem 8.3:} {\em In an $M(QE)_{4}$ spacetime with constant associated scalars $\alpha \neq 0$ and $\beta \neq 0$ satisfying Einstein field equation without cosmological constant the vector field $\xi_{1}$ is a Killing vector field iff the Lie derivative of the energy-momentum tensor with respect to $\xi_{1}$ vanishes.}
\section{Ricci soliton structure in a perfect fluid $M(QE)_{4}$  spacetime} \label{s9}
$~~~~$The idea of Ricci solitons was introduced by Hamilton $\cite{jmh76}$. Ricci solitons also correspond to selfsimilar solutions of Hamilton’s Ricci flow. They are natural generalizations of Einstein metrics and is given by the equation
\begin{equation} \label{9.1}
(L_{V}g)(X,Y)+2S(X,Y)+2\lambda g(X,Y)=0,
\end{equation}
where $\lambda$ is a constant and $V \in \rchi(M)$. There are $3$ categorisations of the Ricci flow according to the values of $\lambda$. Namely,\\\\
\textbf{i) Shrinking:} for $\lambda < 0$,\\
\textbf{ii) Steady:} for $\lambda = 0$,\\
\textbf{ii) Expanding:} for $\lambda > 0$.\\\\
$~~~~$Now for a perfect fluid $M(QE)_{4}$ spacetime using the equations (\ref{8.1}), (\ref{8.4}) and (\ref{8.8}) we get
\begin{equation} \label{9.2}
S(Y,Z)=k\sigma g(Y,Z).
\end{equation}
In the view of equation (\ref{9.2}) the equation (\ref{9.1}) takes the form
\begin{equation} \label{9.3}
(L_{V}g)(Y,Z)=-2(\lambda+k\sigma) g(Y,Z),
\end{equation}
If we consider $\sigma$ as a constant then taking derivative to both the sides of the equation (\ref{9.3}) with respect to $X$ we get
\begin{equation} \label{9.4}
(\nabla_{X}L_{V}g)(Y,Z)=0,
\end{equation}
Now from the commutation formula $\cite{jmh77}$ given by
 \begin{eqnarray} \label{9.5}
&&(L_{V}\nabla_{X}g-\nabla_{X}L_{V}g-\nabla_{[V,X]}g)(Y,Z)\nonumber\\
&&=-g((L_{V}\nabla)(X,Y),Z)-g((L_{V}\nabla)(X,Z),Y),
\end{eqnarray}
the equation (\ref{9.4}) can be written as
\begin{equation} \label{9.6}
g((L_{V}\nabla)(X,Y),Z)+g((L_{V}\nabla)(X,Z),Y)=0.
\end{equation}
A straightforward combinatorial computation yields
\begin{equation} \label{9.7}
(L_{V}\nabla)(Y,Z)=0.
\end{equation}
We know $\cite{jmh77}$ another well known formula
\begin{equation}\label{9.8}
(L_{V}\nabla)(X,Y)=\nabla_{X}\nabla_{Y}V-\nabla_{\nabla_{X}Y}V+R(V,X)Y.
\end{equation}
In view of the formula (\ref{9.8}) the equation (\ref{9.7}) takes the form
\begin{equation} \label{9.9}
\nabla_{Y}\nabla_{Z}V-\nabla_{\nabla_{Y}Z}V+R(V,Y)Z=0.
\end{equation}
Putting $Y=Z=\xi_{1}$ in (\ref{9.9}) we get
\begin{equation} \label{9.10}
\nabla_{\xi_{1}}\nabla_{\xi_{1}}V+R(V,\xi_{1})\xi_{1}=0.
\end{equation}
Hence we get that $V$ is Jacobi along $\xi_{1}$. Thus we obtain the following theorem as:\\\\
\textbf{Theorem 9.1} {\em In a perfect fluid $M(QE)_{4}$ spacetime with constant energy density satisfying Einstein field equation without cosmological constant and admitting a non-trivial(non-Einstein) Ricci soliton, the vector field $V$ is Jacobi along the geodesics of the velocity vector field $\xi_{1}$.}\\\\
Now consider the spacetime is semiconformally flat. Then the equation (\ref{3.3}) can be written as
\begin{eqnarray}\label{9.11}
&&R(X,Y)Z=\alpha\{1+\frac{4m}{3l}\}[g(Y,Z)X-g(X,Z)Y]\nonumber\\
&&~~~~~~~~~~~~~~~~~~~~~+\frac{\beta}{2}[\{A(Y)B(Z)+A(Z)B(Y)\}X\nonumber\\
&&~~~~~~~~~~~~~~~~~~~~~-\{A(X)B(Z)+A(Z)B(X)\}Y\nonumber\\
&&~~~~~~~~~~~~~~~~~~~~~+g(Y,Z)\{A(X)\xi_{2}+B(X)\xi_{1}\}\nonumber\\
&&~~~~~~~~~~~~~~~~~~~~~-g(X,Z)\{A(Y)\xi_{2}+B(Y)\xi_{1}\}].
\end{eqnarray}
Putting $Y=Z=\xi_{1}$ in (\ref{9.11}) implies
\begin{equation}\label{9.12}
R(X,\xi_{1})\xi_{1}=-\alpha\{1+\frac{4m}{3l}\}[X+g(X,\xi_{1})\xi_{1}].
\end{equation}
Taking Lie derivative to both the sides of the equation (\ref{9.12}) we get
\begin{eqnarray}\label{9.13}
&&(L_{V}R)(X,\xi_{1})\xi_{1}+R(X,L_{V}\xi_{1})\xi_{1}+R(X,\xi_{1})L_{V}\xi_{1}\nonumber\\
&&=-\alpha(1+\frac{4m}{3l})[(L_{V}g)(X,\xi_{1})\xi_{1}+g(X,\xi_{1})L_{V}\xi_{1}].
\end{eqnarray}
Taking derivative with respect to $X$ to both the sides of the equation (\ref{9.7}) we get
\begin{equation}\label{9.14}
(\nabla_{X}L_{V}\nabla)(Y,Z)=0.
\end{equation}
The following commutation formula was given by Yano $\cite{jmh77}$,
\begin{equation}\label{9.15}
(L_{V}R)(X,Y)Z=(\nabla_{X}L_{V}\nabla)(Y,Z)-(\nabla_{Y}L_{V}\nabla)(X,Z).
\end{equation}
Equations (\ref{9.14}) and (\ref{9.15}) gives 
\begin{equation}\label{9.16}
(L_{V}R)(X,Y)Z=0.
\end{equation}
Plugging $Y=Z=\xi_{1}$ in (\ref{9.16}) gives
\begin{equation}\label{9.17}
(L_{V}R)(X,\xi_{1})\xi_{1}=0.
\end{equation}
From the equations (\ref{9.3}), (\ref{9.13}), (\ref{9.17}) we get
\begin{eqnarray}\label{9.18}
&&R(X,L_{V}\xi_{1})\xi_{1}+R(X,\xi_{1})L_{V}\xi_{1}\nonumber\\
&&=\alpha(1+\frac{4m}{3l})[2(\lambda+k\sigma)g(X,\xi_{1})\xi_{1}-g(X,\xi_{1})L_{V}\xi_{1}].
\end{eqnarray}
Again putting $Y=Z=\xi_{1}$ in (\ref{9.3}) we get
\begin{equation}\label{9.19}
(L_{V}g)(\xi_{1},\xi_{1})=2(\lambda+k\sigma),
\end{equation}
which implies
\begin{equation}\label{9.20}
g(L_{V}\xi_{1},\xi_{1})=-(\lambda+k\sigma).
\end{equation}
Now taking contraction over $X$ in the equation (\ref{9.18}) and using (\ref{9.20}) we get
\begin{equation}\label{9.21}
-\alpha(1+\frac{4m}{3l})(\lambda+k\sigma)=0.
\end{equation}
If we consider $\alpha(3l+4m)\neq 0$ then equation (\ref{9.21}) will imply $\lambda +k\sigma=0.$ Using the equation (\ref{8.8}) we get
\begin{equation}\label{9.22}
\lambda=kp.
\end{equation}
Hence we derive the following theorem:\\\\
\textbf{Theorem 9.2:} {\em In a semiconformally flat $M(QE)_{4}$ spacetime with constant energy density and with the condition $\alpha(3l+4m) \neq 0,$ satisfying Einstein field equation without cosmological constant and admitting a non-trivial(non-Einstein) Ricci soliton, the following conditions hold:\\
i) the soliton is shrinking if $kp<0$,\\
ii) the soliton is steady if $kp=0$\\
iii) the soliton is expanding if $kp>0$,\\
where $k$ is the gravitational constant and $p$ is the isotropic pressure of the spacetime.}
\section{Example of an $M(QE)_{4}$ spacetime} \label{s10}
$~~~~$Finally, we prove the existence of an $M(QE)_{4}$ spacetime by setting up an example. For this we consider a metric known as Lorentzian metric $g$ on $M^{4}$ given by \\\\
$~~~~~~~~~~~~~$$ds^2 =g_{ij}dx^idx^j= (dx^1)^2+2(x^1)^2(dx^2)^2+3(x^2)^2(dx^3)^2-4(dx^4)^2,$\\\\
where $i,j=1,2,3,4.$ Here we obtain the non-vanishing components of Christoffel symbols, curvature tensors and Ricci tensors which as follows:\\\\
\begin{equation}\label{10.1}
\renewcommand{\arraystretch}{1.5}
\begin{array}{r@{\;}l}
\Gamma_{33}^2=-\frac{3}{2}\frac{(x^2)}{(x^1)^2},~\Gamma_{22}^1=-2(x^1),~\Gamma_{23}^3=\frac{1}{(x^2)},~\Gamma_{12}^2=\frac{1}{(x^1)},\\
\end{array}
\end{equation}
\begin{eqnarray}\label{10.2}
\renewcommand{\arraystretch}{1.5}
\begin{array}{r@{\;}l}
R_{1332}=-3\frac{(x^2)}{(x^1)},~S_{12}=-\frac{1}{(x^1)(x^2)}.\\
\end{array}
\end{eqnarray}
Now we wish to prove that this is an $M(QE)_{4}$ manifold.\\\\
$~~~~$We consider $\alpha, \beta$ as associated scalars and we take these scalars as follows:
\begin{equation}\label{10.3}
\alpha=0,~\beta=-\frac{\sqrt{3}}{2(x^1)^2(x^2)}\neq 0,
\end{equation}
and the associated $1$-forms are as follows:\\\\
\begin{math}
~~~~~~~~~~~~A_{i}(x)=\left \{\begin{array}{cccl}
{0} & \mbox{for} & i=1,3,4 \\
{\sqrt{2}(x^1)} & \mbox{for} & i=2 
\end{array}
\right.
\end{math}~~~~~~;~~~~~~~\begin{math}B_{i}(x)=\left
\{\begin{array}{ccl}
{\sqrt{\frac{2}{3}}} & \mbox{for} & i=1 \\
{(x^2)} & \mbox{for} & i=3 \\
{0} & \mbox{for} & i=2,4.
\end{array}
\right.
\end{math}\\\\\\
By using the equations (\ref{10.2}) and (\ref{10.3}) and using the values of $A_{i},~B_{i}$ as described above, we get the following equation:
\begin{equation}\label{10.4}
S_{12}=\alpha g_{12}+\beta [A_{1}B_{2}+B_{1}A_{2}],
\end{equation}
which matches with our desired equation described by (\ref{2.5}). Except $S_{12}$ any other $S_{ij}$ is $0~\forall~i,j=1,2,3,4$. On the other hand the following equations can also be verified very easily:
\begin{equation}\label{10.5}
g^{ij}A_{i}A_{j}=1,~g^{ij}B_{i}B_{j}=1,~g^{ij}A_{i}B_{j}=0~\forall~i,j=1,2,3,4,
\end{equation}
i.e. the $1$-forms $A$ and $B$ are both unit and they are also orthogonal to each other. Hence we obtain\\\\
$~~~~~~~~~~~~~~~~~~~~~~~~$$S_{ij}= \alpha g_{ij}+\beta[A_{i}B_{j}+B_{i}A_{j}],~\forall~i,j= 1,2,3,4,$\\\\
where $A$ and $B$ are two nonzero unit $1$-forms orthogonal to each other and $\beta \neq 0$. Thus $(M^{4}, g)$ is a mixed quasi-Einstein spacetime.\\\\
\textbf{Example 10.1.} \emph{Let $(M^{4}, g)$ be a Lorentzian manifold furnished with the Lorentzian metric $g$ given by\\\\
$~~~~~~~~~~~~~$$ds^2 =g_{ij}dx^idx^j= (dx^1)^2+2(x^1)^2(dx^2)^2+3(x^2)^2(dx^3)^2-4(dx^4)^2,$\\\\
where $i,j=1,2,3,4.$ Then $(M^{4}, g)$ is an $M(QE)_{4}$ spacetime.}\\\\
\textbf{Conclusion:} {\em In the last few decades there have been rigorous inspections in the mathematical theory to study the physical properties of this universe. Many mathematicians/mathematical physicists have worked tirelessly to generalize and study the nature of this universe mathematically. Our present work is a new note to add something more to that knowledge. Here we have worked on both mathematical and physical significance of a mixed quasi-Einstein spacetime. In the first part of our paper we have worked on the geometry of a mixed quasi-Einstein spacetime and in the last part we have drawn some physical applications of the same spacetime also. Since it works as a bridge between an Einstein spacetime and a generalised quasi-Einstein spacetime, which we need to comprehend the second phase of the evolution of the universe, thus the present work is relevant and contemporary.}

\end{document}